\documentclass[12pt]{article}
\usepackage{amsfonts,amssymb}
\usepackage{latexsym}
\usepackage{theorem}
\IfFileExists{mathrsfs.sty}{\usepackage{mathrsfs}\let\mathcal\mathscr}{}
\newcommand{\linespacing}{1.25}
\renewcommand{\baselinestretch}{\linespacing}
\newtheorem{theorem}{Theorem}[section]
\newtheorem{lemma}[theorem]{Lemma}
\newtheorem{proposition}[theorem]{Proposition}
\newtheorem{corollary}[theorem]{Corollary}

{\theorembodyfont{\normalfont\rmfamily}
\newtheorem{definition}[theorem]{Definition}
\newtheorem{remark}[theorem]{Remark}
}
\makeatletter
\def\ack{\vspace{.5\baselineskip}\noindent{\theorem@headerfont
Acknowledgement}\ \ }
   {\renewcommand{\theththm}{\ref{#1}}\begin{ththm}}
   {\end{ththm}}
\newtheorem{ththm}{Theorem}
\newenvironment{proof}[1][]%
{\def\proof@temp{#1}\par\noindent
\textsc{Proof}\ifx\proof@temp\@empty\else\
({\proof@temp})\fi\hspace{1em}}
{\hphantom{xxx}\hfill~ {$\Box$}\par\vspace{.4\baselineskip}}
\def\operatorname#1{\mathop{\operator@font #1}\nolimits}%
\makeatother
\newcommand{\map}[1]{\stackrel{#1}{\longrightarrow}}
\newcommand{\PO}{P}
\newcommand{\CiM}{C^\infty(M)}
\newcommand{\bbnu}{[\![\nu]\!]}
\newcommand{\Exp}{\operatorname{Exp}}
\newcommand{\C}{\mathbb{C}}
\newcommand{\R}{\mathbb{R}}
\newcommand{\K}{\mathbb{K}}
\newcommand{\W}{\mathcal{W}}
\let\wt\widetilde

\let\ol\overline
\newcommand{\cob}{\mathop{\mathfrak{d}}}
\newcommand{\D}{\mathcal{D}}
\renewcommand{\L}{\mathcal{L}}
\newcommand{\BR}{\{\, ,\,\}}

\newcommand{\g}{\mathfrak{g}}

\newcommand{\ad}{\operatorname{ad}}

\newcommand{\Der}{\operatorname{Der}}

\newcommand{\Id}{\operatorname{Id}}
\newcommand{\Inn}{\operatorname{Inn}}

\newcommand{\half}{{\textstyle{\frac12}}}
\def\cyclic{\mathop{\kern0.9ex{{+}
\kern-2.2ex\raise-.28ex\hbox{\Large\hbox
{$\circlearrowright$}}}}\limits}

\def\ftnote#1{\def\footnotemark{}\footnote{#1}\setcounter{footnote}{0}}

\begin{document}

\title{Natural star products on symplectic manifolds \\and
quantum moment maps
\ftnote{This research was partially supported by an
Action de Recherche Concert\'ee de la Communaut\'e
fran\c{c}aise de Belgique and by the Belgian FNRS.}\\\ }

\author{Simone Gutt${}^{1,2}$\ftnote{\kern-4.5pt${}^{1}$Universit\'e Libre 
de Bruxelles, Campus Plaine, CP 218, BE-1050~Brussels, Belgium}
\ftnote{\kern-4.5pt${}^{2}$Universit\'e de Metz, Ile du Saulcy,
F-57045~Metz Cedex 01, France}
\\[20pt]
\& \\[20pt]
John Rawnsley${}^{1,3}$\ftnote{\kern-4.5pt${}^{3}$Mathematics Institute,
University of Warwick, Coventry CV4 7AL, United Kingdom}\ftnote{Email: 
\texttt{sgutt@ulb.ac.be}, \texttt{j.rawnsley@warwick.ac.uk}}
\\[10pt]\ }

\date{\textit{It is a pleasure to dedicate this paper to Alan Weinstein\\ 
on the occasion of  his sixtieth birthday.}}
\setcounter{page}{0}

\renewcommand{\baselinestretch}{1}

\maketitle

\thispagestyle{empty}

\begin{abstract}
We define a natural class of star products: those which are given by a
series of bidifferential operators which at order $k$ in the deformation
parameter have at most $k$ derivatives in each argument. We show that
any such star product on a symplectic manifold defines a unique
symplectic connection. We parametrise such star products, study their
invariance and give necessary and sufficient conditions for them to
yield  a quantum moment map.

We show that Kravchenko's sufficient condition \cite{bib:Kravchenko} for
a moment map for a Fedosov star product is also necessary.
\end{abstract}

\newpage

\renewcommand{\baselinestretch}{\linespacing}

\section{Introduction}\label{sect:intro}
The relation between a star product on a symplectic manifold and
a symplectic connection on that manifold appears in many contexts.
In particular, when one studies properties of
invariance of star products, results are much easier when there is
an invariant connection.
We show here that there is a natural class of star products which
define a unique symplectic connection.
We study the invariance of such products and the conditions for
them to have a moment map.

\ack
The results on moment maps were presented
in Bayrischzell in April 2002 (and have been extended
by M\"uller and Neumaier \cite{bib:MulNeu}), those
on natural star products were presented in Zurich in January 2003;
the first author thanks the organisers for their
hospitality. The second author thanks the Belgian FNRS for its support
when visiting Brussels for this project.

\section{Natural star products}
Let $(M,\PO)$ be a Poisson manifold. Let $\CiM$ be the space of
$\K$-valued smooth functions on $M$, where $\K$ is $\R$ or $\C$. The
space of Hochschild $(p+1)$-cochains on $\CiM$ with values in $\CiM$
which are given by differential operators in each argument
(\textbf{differential $(p+1)$-cochains}) will be denoted by $\D^p$.
Those which are differential operators of order at most $k$ in each
argument will be denoted by $\D_k^p$. $\cob \,\colon \D^p \to \D^{p+1}$
denotes the Hochschild coboundary operator, and $[A,B]$ denotes the
Gerstenhaber bracket of cochains which we also write $\ad A\
B$. If $m(u,v)=uv$ denotes the multiplication of functions and $br(u,v)
=
\{u,v\}=\PO(du,dv)$ the Poisson bracket then $m \in \D^1_0$ and $br \in
\D^1_1$. Finally we have the relation $\cob A = - \ad m\ A$ and
$[\D^p_r,\D^q_s] \subset \D^{p+q}_{r+s-1}$.

\begin{definition}\label{def:natstar}
A \textbf{natural star product} on $(M,\PO)$ is a bilinear map
\[
\CiM\times \CiM \to \CiM\bbnu, \qquad (u,v) \mapsto u*v : =
\sum_{r\ge0} \nu^rC_r(u,v),
\]
which defines a formally
associative product on $\CiM\bbnu$
   when the map is extended $\K\bbnu$-linearly
(i.e.\ $(u*v)*w = u*(v*w))$ and such that
\begin{itemize}
\item $C_0=m$;
\item  the skewsymmetric part of $C_1$ is the Poisson bracket:
   $C_1(u,v)-C_1(v,u) = 2\{u,v\}$;
\item  $1*u = u*1 = u$;
\item  each $C_r$ is in $\D_r^1$ (i.e.\ is a  \textbf{bidifferential
operator} on $M$
    \textbf{of order at most} $r$  \textbf{in each argument}).
\end{itemize}
\end{definition}

\begin{remark}
It would be enough for some purposes to assume that $C_1$ be of order 1,
and $C_2$ of order at most 2 in each argument; c.f.{} \cite{bib:Xu}
where natural is defined this way, except that $C_1$ is the Poisson
bracket. The main results still hold with the appropriate modifications.
\end{remark}

\begin{proposition}\label{prop:equiv}
Two natural star products $*$ and $*'$ on $(M,P)$ are equivalent if and
only if there is a series
\[
E = \sum_{r = 1}^\infty \nu^r E_r
\]
where the $E_r$ are differential operators of order at most $r+1$,
such that
\begin{equation}\label{intro:nateq}
f*'g = \Exp E\,((\Exp -E)\, f * (\Exp -E)\,g)),
\end{equation}
where $\Exp{}$ denotes the exponential series.
\end{proposition}

\begin{proof}
We have, for any $E \in \D^0$, $C \in \D^1$
\[
(\ad E \, C)(u,v)=E(C(u,v))-C(Eu,v)-C(u,Ev).
\]
and so
\[
(\,\Exp (\ad E)\, C)\,(u,v)= \Exp E\,(C\,((\Exp -E)\, u,(\Exp -E)\, v)).
\]

Consider now a natural star product $*=\sum_{r\ge0} \nu^rC_r$ and a
series $E = \sum_{r = 1}^\infty \nu^r E_r$ where the $E_r$ are
differential operators of order at most $r+1$. From the observation
above and $[\D^0_r,\D^1_s] \subset \D^1_{r+s-1}$ we have
\[
f*'g = \Exp E\,(\Exp -E f* \Exp -E g)= (\Exp (\ad E)\, *)(u,v)
\]
is a natural star product on $M$.

Reciprocally, if $*$ and $*'$ are two natural star products which are
equivalent, we shall show that the equivalence is necessarily given by
$f*'g = \Exp E\,(\Exp -E \, f* \Exp -E \, g)= (\Exp (\ad E)\, *)(u,v)$
where  $E = \sum_{r = 1}^\infty \nu^r E_r$ with $E_r\in\D_{r+1}^0$. The
fact that $*$ and $*'$ are equivalent can be written $f*'g = \Exp
E\,(\Exp -E f* \Exp -E g)$ where $E = \sum_{r = 1}^\infty \nu^r E_r$ is
a formal series of linear maps from $\CiM$ to $\CiM$. We shall prove by
induction that all the $E_r$ are in $D_{r+1}^0$. Assume the $E_r$ are in
$\D_{r+1}^0$ for $r\le k$. Define $E^{(k)} = \sum_{r = 1}^k \nu^r E_r$
and $*^{''}= (\Exp (\ad E^{(k)})\, *)=:\sum_{r\ge0} \nu^rC^{''}_r$. It
coincides with $*'$ up to and including order $k$, and it is a natural star
product on $M$. The difference $C^{''}_{k+1}-C'_{k+1}$ is a
bidifferential operator which is a Hochschild $2$-cocycle whose
skewsymmetric part vanishes because it is given by $\cob E_{k+1}$; hence
it is the   coboundary of a differential operator; since this
bidifferential operator is of order at most $k+1$ in each argument, it
is the coboundary $=\cob E'_{k+1}$ of an element $ E'_{k+1} \in
\D_{k+2}^0$. Thus $E_{k+1}-E'_{k+1}$ is a $1$-cocycle, i.e.\ a vector
field. This proves that $E_{k+1}$ is in $\D_{k+2}^0$.
   \end{proof}

\begin{remark}
The space of formal differential operators considered in Proposition
\ref{prop:equiv}, $E = \sum_{r = 1}^\infty \nu^r E_r$ where the $E_r$
are differential operators of order at most $r+1$, is a pronilpotent Lie
algebra.
\end{remark}

\begin{remark}
All the explicit constructions of star products are natural:
\begin{itemize}
\item star products on cotangent bundles
\cite{bib:BordNeuWald,bib:CahGutt,
bib:DeWildeLecomte};
\item  star products given by Fedosov's construction
\cite{bib:Fedosov1,bib:Fedosov2};
\item  star products on the dual of a Lie algebra or coadjoint orbits
\cite{bib:ArnalCortet,bib:Gutt};
\item star products with separation of variables on a Kaehler manifold
\cite{bib:Karabegov};
\item star products given by Kontsevich's construction
\cite{bib:Kontsevich}.
\end{itemize}
\end{remark}

\section{Connections}\label{sect:conn}

The link between the notion of star product on a symplectic manifold and
symplectic connections already appears in the seminal paper of Bayen,
Flato, Fronsdal, Lichnerowicz and Sternheimer \cite{bib:Bayen}, and was
further developed by Lichnerowicz \cite{bib:Lichne} who showed that any
so called Vey star product (i.e.\ a star product defined by
bidifferential operators whose principal symbols at each order coincide
with those of the Moyal star product) determines a unique symplectic
connection. Fedosov gave a construction of Vey star products starting
{}from a symplectic connection and a series of closed two forms on the
manifold. It was shown that any star product is equivalent to a Fedosov
star product. Nevertheless, many star products which appear in natural
contexts (cotangent bundles, Kaehler manifolds\dots) are not Vey star
products (but are natural in the sense defined above). The aim of this
section is to generalise the result of Lichnerowicz and to show that
on any symplectic manifold, a natural star product determines a unique
symplectic connection.

Given any torsion free linear connection $\nabla$ on $(M,P)$,
the  term of order $1$ of a natural star product can be written
\[
C_1=\BR-\cob E_1=\BR + (\ad E_1)\,m \qquad \mbox{\rm{where }}
\, E_1\in \D^0_2
\]
and the term of order $2$ can be written in a chart
\begin{eqnarray*}
C_2(u,v)&=&\half ((\ad E_1)^2\,m)(u,v) +((\ad E_1)\,\BR)(u,v)\\
&&\quad\mbox{}+\half \PO^{ij}\PO^{i'j'}\,\nabla^{2}_{ii'} u\,
\nabla^{2}_{jj'}v\\
&&\qquad\mbox{}+ {\textstyle{\frac 16}}( \PO^{rk}\nabla_r\PO^{jl}+
   \PO^{rl}\nabla_r\PO^{jk})(\nabla^{2}_{kl} u\, \nabla_{j}v
   +\,\nabla_{j} u\, \nabla^{2}_{kl}v)\\
&&\qquad\quad-\cob E_2(u,v) +c_2(u,v)
\end{eqnarray*}
where $ E_2\in \D^0_3$ and where $c_2\in \D^1_1$ is skewsymmetric.

Remark that $E_1$ is not uniquely  defined; two choices differ by
an element $X \in \D^0_1$. Observe that the first lines in the
definition of $C_2$ for two such different choices only differ
by an element in $\D^1_1$. Indeed
\begin{eqnarray*}
\half \left(\ad (E_1+X)^2\,m\right)(u,v) +
\left((\ad E_1+X)\,\BR\right)(u,v)
&=&\\
   &&\kern-2in\half ((\ad E_1)^2\,m)(u,v) +((\ad E_1)\,\BR)(u,v)\\
  && \kern-2in\quad\mbox{}+\half ((\ad X)\circ(\ad E_1)\,m)(u,v) +
  ((\ad X)\,\BR)(u,v)\
   \end{eqnarray*}
and $(\ad E_1)\,m$, $(\ad X)\,\BR$ are in $\D_1^1$,
so also is $((\ad X)\circ(\ad E_1)\,m)$.

Changing the torsion free linear connection gives a modification
of the terms of the second line of $C_2$; writing
$\nabla'=\nabla + S$, this modification involves
terms of order $2$ in one argument and $1$ in the other given by
\begin{eqnarray*}
\left(-\half \PO^{rk}\PO^{^sl}S_{rs}^j
+ {\textstyle{\frac 13}}( \PO^{rk}S_{rs}^j\PO^{sl}
                              +\PO^{rk}S_{rs}^l\PO^{js})\right)
   (\nabla^{2}_{kl} u\, \nabla_{j}v
   +\,\nabla_{j} u\, \nabla^{2}_{kl}v)
   &=&\\
   &&\kern-2.5in\mbox{}-\cyclic_{jkl}{\textstyle{\frac 16}}
   \PO^{rk}\PO^{sl}S_{rs}^j
   (\nabla^{2}_{kl} u\, \nabla_{j}v
   +\,\nabla_{j} u\, \nabla^{2}_{kl}v)\\
\end{eqnarray*}
as well as terms of order $1$ in each argument, where $\cyclic$ denotes
a cyclic sum over the indicated variables.

Notice that the terms above coincide with the terms of
the same order in the coboundary of the operator
$E'={\textstyle{\frac 16}}\cyclic_{jkl} \PO^{rk}\PO^{sl}S_{rs}^j
   \nabla^{3}_{jkl}$.

If the Poisson tensor is invertible (i.e.\ in the symplectic
situation), the symbol of any differential operator of order
$3$ can be written in this form $E'$, hence we have:

\begin{proposition}
A  star product $*=\sum_{r\ge0}\nu^rC_r$ on a symplectic manifold
$(M,\omega)$, so that $C_1$ is a bidifferential operator of order
 1 in each argument and $C_2$ of order at most 2 in each argument,
 determines a unique symplectic connection $\nabla$ such
that
\begin{equation}\label{intro:defconn}
C_1=\BR- \cob E_1
\qquad C_2=\half (\ad E_1)^2\,m + ((\ad E_1)\,\BR)+ \half\PO^{2}(\nabla^{2}
\cdot,\nabla^{2}\cdot)
+A_2
\end{equation}
where $A_2 \in \D^1_1$ and
$\PO^{2}(\nabla^{2}u,\nabla^{2}v)$ denotes the bidifferential
operator which is given by
$\PO^{ij}\PO^{i'j'}\,\nabla^{2}_{ii'} u\, \nabla^{2}_{jj'}v$
in a chart.
\end{proposition}

\begin{remark}
This shows, in particular, that any natural star product
$*=\sum_{r\ge0}\nu^rC_r$ on a symplectic manifold
$(M,\omega)$ determines a unique symplectic connection $\nabla$.
\end{remark}

\section{Symmetry and Invariance}\label{sect:symm}

Symmetries in quantum theories are automorphisms of the algebra
of observables. Thus we define a symmetry $\sigma$ of a star product
$* = \sum_r \nu^rC_r$ as an automorphism of the $\K\bbnu$-algebra
$\CiM\bbnu$
with multiplication given by $*$:
\[
\sigma(u*v) = \sigma(u) * \sigma(v), \qquad \sigma(1) =1,
\]
where $\sigma$, being determined by what it does on $\CiM$, 
will be a formal series
\[
\sigma(u) = \sum_{r\ge0} \nu^r \sigma_r(u)
\]
of linear maps $\sigma_r \colon \CiM \to \CiM$. In terms of the
components of $\sigma$, the conditions to be an automorphism are
\begin{eqnarray}
&\circ&\label{eq:auto}
\sum_{r+s=k} \sigma_r(C_s(u,v))
= \sum_{r+s+t=k} C_r(\sigma_s(u), \sigma_t(v)),\quad k\ge0,
 u,v \in \CiM;\\
&\circ&
\hbox{$\sigma_0$ is invertible;}\nonumber\\
&\circ&\label{eq:auto1}
\sigma_0(1)=1,\qquad \sigma_r(1) = 0, \quad r\ge1.
\end{eqnarray}

\begin{lemma}
If $*$ is a star product on a Poisson manifold $(M,P)$ and $\sigma$ is
an automorphism of $*$ then it can be written $\sigma(u) = T(u \circ
\tau)$ where $\tau$ is a Poisson diffeomorphism of $(M,P)$ and $T = \Id
+ \sum_{r\ge1} \nu^r T_r$ is a formal series of linear maps. If $*$ is
differential, then the $T_r$ are differential operators; if $*$ is
natural, then $T= \Exp E$ with $E = \sum_{r\ge1} \nu^r E_r$ and $E_r$
is a
differential operator of order at most $r+1$.
\end{lemma}

\begin{proof}
Taking $k=0$ in (\ref{eq:auto}), we have $\sigma_0(uv) = \sigma_0(u)
\sigma_0(v)$ so $\sigma_0$ is an automorphism of $\CiM$ and hence is
composition with a diffeomorphism $\tau$ of $M$. Taking $k=1$ in
(\ref{eq:auto}), and antisymmetrising in $u$ and $v$ we have
$\sigma_0(\{u,v\}) = \{\sigma_0(u),\sigma_0(v)\}$, and so $\tau$ is a
Poisson map.

Set $T(u) = \sigma(u \circ \tau^{-1})$ so $T = \sigma \circ
\sigma_0^{-1}$ which has the stated form. Define a new star product
by $u*'v = \sigma_0( \sigma_0^{-1}(u) * \sigma_0^{-1}(v))$ then $*'$
is differential (resp.~natural) if $*$ is. On the other hand
$\sigma(u*v) = \sigma(u) * \sigma (v)$ implies $ T(\sigma_0(u*v)) =
T(\sigma_0(u)) * T(\sigma_0 (v))$ and hence that $T(u*'v) = T(u) *
T(v)$.
Thus $T$ is an equivalence between $*'$ and $*$. The result now follows
{}from Theorem 2.22 of \cite{bib:GuttRaw} and
Proposition~\ref{prop:equiv}.
\end{proof}

If $\sigma_t$ is a one-parameter group of symmetries of the star product
$*$, then its generator $D$ will be a derivation of $*$. Denote the Lie
algebra of derivations of $*$ by $\Der(M,*)$. Moreover by
differentiating
the statement of the Lemma, $D = \sum_{r \ge 0} \nu^r D_r $ with $D_0 =
X$, a Poisson vector field ($\L_X P = 0$), and if $*$ is natural then
each $D_r$ for $r \ge 1$ is a differential operator of order at most
$r+1$. In the case we have a Lie group $G$ acting by symmetries, then
there will be an induced action of $G$ on $M$ and
the infinitesimal automorphisms will give a homomorphism of Lie algebras $D
\colon \g \to \Der (M,*)$ from its Lie algebra $\g$ into the derivations
of the star product of the above form. For each $\xi \in \g$, $D_\xi =
\wt{\xi} + \sum_{r\ge1} \nu^r D_\xi^r$ where $\wt{\xi}$ is the vector
field generating the induced action of $\exp t \xi$ on $M$.

\begin{definition}\label{def:inv}
A star product $*=m+\sum_{r\ge1}\nu^rC_r$ on a Poisson manifold $(M,P)$
is said to be \textbf{invariant} under a diffeomorphism $\tau$ of $M$ if
$u \mapsto u \circ \tau$ is a symmetry of $*$.
\end{definition}

Observe that if $*$ is $\tau$-invariant then $\tau$ preserves each
cochain $C_r$ and hence the Poisson bracket. Invariance is a much
stronger condition than being the leading term of a symmetry.

\section{Quantum moment maps}\label{sect:qmm}

A derivation $D \in \Der(M,*)$ is said to be \textbf{essentially inner}
or
\textbf{Hamiltonian} if $D = \frac1\nu \ad_* u$ for some $ u \in
\CiM\bbnu$. We denote by $\Inn(M,*)$ the essentially inner derivations
of
$*$. It is a linear subspace of $\Der(M,*)$ and is the quantum analogue
of
the Hamiltonian vector fields.

By analogy with the classical case, we call an action of a Lie group
\textbf{almost $*$-Hamiltonian} if each $D_\xi$ is essentially inner, and
call a linear choice of functions $u_\xi$ satisfying
\[
D_\xi = \textstyle{\frac1\nu} \ad_* u_\xi, \qquad \xi \in \g
\]
a (quantum) Hamiltonian.
We say the action is \textbf{$*$-Hamiltonian} if $u_\xi$ can be chosen to
make $\xi \mapsto u_\xi \colon \g \to \CiM\bbnu$ a homomorphism of Lie
algebras. When $D_\xi=\wt{\xi}$, this map is called \textbf{a quantum moment
map} \cite{bib:Xu}.

Considering the map $ a \colon \CiM\bbnu \to \Der(M,*)$ given by
\[
a(u) (v)= \textstyle{\frac1\nu} \ad_* u(v)
= \textstyle{\frac1\nu} (u*v - v*u).
\]
and defining a bracket on $\CiM\bbnu$ by
\[
[u,v]_* = \textstyle{\frac1\nu} (u*v - v*u)
\]
then, by associativity of the star product, $a$ is a homomorphism of Lie
algebras whose image is $\Inn(M,*)$. Since $D\circ a(u) - a(u) \circ
D = a(Du)$, $\Inn(M,*)$ is an ideal in $\Der(M,*)$ and so there is an
induced Lie bracket on the quotient $\Der(M,*)/\Inn(M,*)$.

\begin{lemma}(\cite{bib:BertCahGutt})
If $*$ is a star product on a symplectic manifold $(M,\omega)$ then
the space of derivations modulo inner derivations,
$\Der(M,*)/\Inn(M,*)$, can be identified with $H^1(M,\R)\bbnu$ and the
induced bracket is zero.
\end{lemma}

\begin{proof}
The first part is well known; let us recall that
 locally any derivation $D \in \Der(M,*)$ is
inner, and that the ambiguity in the choice of
a corresponding function $u$ is locally constant so that the exact
1-forms $du$ agree on overlaps and yield a globally defined (formal)
closed 1-form $\alpha_D$. The map $\Der(M,*) \to
Z^1(M)\bbnu$ defined by $D \mapsto \alpha_D$ if $D|_U =
\textstyle{\frac1\nu} \ad_* u$ and $\alpha_D|_U = du$ is a linear
isomorphism with the space of formal series of closed 1-forms and
maps essentially inner
derivations to exact 1-forms inducing a bijection $\Der(M,*)/\Inn(M,*)
\to Z^1(M)\bbnu/d(\CiM\bbnu) = H^1(M,\R)\bbnu$.

Let $D_1$ and $D_2$ be two derivations of $*$. $(D_1 \circ D_2 - D_2
\circ D_1)|_U = a([u_1,u_2]_*)$. But $[u_1,u_2]_*$ does not change if we
add a local constant to either function, so is the restriction to $U$ of
a globally defined function which depends only on $D_1$ and $D_2$. We
denote this function by $b(D_1,D_2)$ and have the identity
\[
D_1 \circ D_2 - D_2 \circ D_1 = \textstyle{\frac1\nu} \ad_* b(D_1,D_2).
\]
This shows that $[\Der(M,*),\Der(M,*)] \subset \Inn(M,*)$ and hence that
the induced bracket on $H^1(M,\R)\bbnu$ is zero.
\end{proof}

The kernel of $a$ consists of the locally constant formal functions
$H^0(M,\R)\bbnu$ and hence:

\begin{remark}
If $*$ is a differential star product on a symplectic manifold
$(M,\omega)$ then there is an exact sequence of Lie algebras
\[
0 \to H^0(M,\R)\bbnu \hookrightarrow \CiM\bbnu \map{a} \Der(M,*) \map{c}
H^1(M,\R)\bbnu \to 0
\]
where $c(D) = [\alpha_D]$.
\end{remark}

\begin{corollary}(see also \cite{bib:Xu})
Let $G$ be a Lie group of symmetries of a star product $*$ on
$(M,\omega)$ and $d\sigma \colon \g \to \Der(M,*)$ the induced
infinitesimal action. If $H^1(M,\R) = 0$ or $[\g,\g] \subset \g$ then
the action is almost $*$-Hamiltonian.
\end{corollary}

Indeed,
by definition, the action is almost $*$-Hamiltonian if $d\sigma(\g) \subset
\Inn(M,*)$. This is the case under either of the two conditions.

\section{Moment Maps for a Fedosov Star Product}

In this section we examine the necessary and sufficient conditions for a
Fedosov star product to have a moment map. In order to do this it is
necessary to examine the Fedosov construction in detail, which we
therefore repeat here.

Having chosen a series of closed 2-forms $\Omega \in \nu
\Lambda^2(M)\bbnu$ and a symplectic connection $\nabla$ on a symplectic
manifold $(M,\omega)$, we consider the Fedosov star product associated
to these data. This star product is obtained by identifying $\CiM\bbnu$
with the space of flat sections of the Weyl bundle $\W$ (which is bundle
of associative algebras) endowed with a flat connection (the Fedosov
connection).

Sections of the Weyl bundle have the form of formal series
\[
a(x,y,\nu) = \sum_{2k+l\ge0} \nu^k a_{k,i_1,\ldots,i_l}(x) y^{i_1}\cdots
y^{i_l}
\]
where the coefficients $a_{k,i_1,\ldots,i_l}$ are symmetric covariant
tensor fields on $M$; $2k+l$ is the degree in $\W$ of the corresponding
homogeneous component.

The product of two sections taken pointwise makes the space of sections
into an algebra, and in terms of the above representation of sections
the multiplication has the form
\[
(a\circ b)(x,y,\nu) = \left.\left(\Exp
\left(\textstyle{\frac\nu2}\Lambda^{ij}
\frac{\partial}{\partial y^i} \frac{\partial}{\partial z^j} \right)
a(x,y,\nu)b(x,z,\nu)\right)\right|_{y = z},
\]
with $\Lambda^{ij} \omega_{jk} = \delta^i_k$ (thus $\{ f,g\} =
\Lambda^{ij} \partial_i f\partial_j g$).

If we introduce the $\W$-valued 1-form $\ol{\Gamma}$ given by
\[
\ol{\Gamma} = \half \omega_{ki} \Gamma^k_{rj} y^i y^j dx^r
\]
then the connection in $\W$ is given by
\[
\partial a = da - \textstyle{\frac1\nu}[\ol{\Gamma},a].
\]

More generally one looks at forms with values in the Weyl bundle, and
locally sections of  $\W\otimes \Lambda^q$ have the form
\[
\sum_{2k+p\ge0} \nu^k a_{k,i_1,\ldots,i_l,j_1,\ldots,j_q}(x)
y^{i_1}\dots
y^{i_p} \,dx^{j_1}\wedge\dots\wedge dx^{j_q}
\]
where the coefficients are again covariant tensors, symmetric in
$i_1,\ldots,i_p$ and anti-symmetric in $j_1,\ldots,j_q$. Such sections
can be multiplied using the product in $\W$ and simultaneously exterior
multiplication $a\otimes\omega \circ b\otimes \omega^\prime = (a\circ
b)\otimes(\omega \wedge \omega^\prime)$. The space of $\W$-valued forms
$\Gamma(\W\otimes\Lambda^*)$ is then a graded Lie algebra with respect
to the bracket
\[
[s,s^\prime] = s \circ s^\prime - (-1)^{q_1q_2} s^\prime \circ s.
\]

As usual, the connection $\partial$ in $\W$ extends to a covariant
exterior derivative on all of $\Gamma(\W\otimes\Lambda^*)$, also denoted
by $\partial$ by using the Leibnitz rule:
\[
\partial(a\otimes\omega) = \partial (a) \wedge \omega + a \otimes
d\omega.
\]
The curvature of $\partial$ is then given by $\partial\circ\partial$
which is a 2-form with values in $\hbox{\rm End}(\W)$. In this case
it admits a simple expression in terms of the curvature $R$
   of the symplectic connection $\nabla$:
\[
\partial\circ\partial a = \textstyle{\frac1\nu}[\ol{R},a]
\]
where
\begin{equation}\label{R}
\ol{R} = \textstyle{\frac14} \omega_{rl} R^l_{ijk}y^ry^k\,dx^i\wedge
dx^j.
\end{equation}

Define
\[
\delta(a) =  dx^k\wedge \frac{\partial a}{\partial y^k}.
\]
Note that $\delta$ can be written in terms of the algebra structure by
\[
\delta(a) = \textstyle{\frac1\nu}\left[- \omega_{ij}y^idx^j,a\right].
\]

With these preliminaries we construct a connection $D$ on $\W$ of the
form
\begin{equation}\label{D}
Da = \partial a - \delta(a) - \textstyle{\frac1\nu}[r,a]
\end{equation}
which is flat: $D\circ D = 0$. Since
\[
D\circ D a = \textstyle{\frac1\nu}\left[\ol{R} +\delta r-\partial r  +
\textstyle{\frac\nu2}[r,r], a \right]
\]
one takes the solution $r$ so that
\[
\ol{R} +\delta r-\partial r  + \textstyle{\frac1{2\nu}}[r,r] = \Omega
\]
and hence so that
\begin{equation}\label{Dr}
Dr = \ol{R}-\Omega-\textstyle{\frac1{2\nu}}[r,r].
\end{equation}

This solution is given by the iterative process
\[
r = \delta^{-1} \left(-\ol{R} +\partial r-\textstyle{\frac1{2\nu}}[r,r]+
\Omega\right)
\]
with $r$ a $1$-form with values in elements of $\W$ of degree at least
$2$. In the formula above $ \delta^{-1}a_{pq} = \frac1{p+q} y^k
i(\frac{\partial}{\partial x^k})a_{pq}$ if $p+q>0$ and $\delta^*a_{pq} =
0$ if  $p+q = 0$ where $a_{pq}$ denotes the terms in $a$ corresponding
to a $q$-form with $p$ $y$'s. Remark then that $(\delta^{-1}\circ
\delta+\delta\circ\delta^{-1})a = a-a_{00}$.

One has, for any smooth vector field $X$ on $M$:
\[
\delta\circ i(X)+ i(X)\circ \delta = \textstyle{\frac1\nu}\ad_*(
\omega_{ij}X^iy^j)
\]
\[
\ad_* r\circ i(X)+ i(X)\circ \ad_* r = \ad_* (i(X)r)
\]
and
\[
\partial\circ i(X)+ i(X)\circ \partial = \L_X -
(\nabla_{i}X)^j y^i \partial_{y^j}
\]
which can be rewritten as
\[
\partial\circ i(X)+ i(X)\circ \partial = \L_X
+\textstyle{\frac1\nu}\ad_*\left(-\half (\nabla_{i}(i(X)
\omega))_{j}y^iy^j
\right)
+\half  (di(X)\omega)_{ip}y^i\Lambda^{jp}\partial_{y^j}.
\]
This gives the generalised Cartan formula (which coincides with the one
given by Neumaier \cite{bib:Neumaier})
\begin{eqnarray}
\L_X
& = &D\circ i(X)+ i(X)\circ D+\textstyle{\frac1\nu}\ad_*(
\omega_{ij}X^iy^j)
+\textstyle{\frac1\nu} \ad_* (i(X)r)\\
&&\qquad\mbox{}+\textstyle{\frac1\nu} \ad_*\left(\half (\nabla_{i}(i(X)
\omega))_{j}y^iy^j\right) - \half  (di(X)\omega)_{ip} y^i\Lambda^{jp}
\partial_{y^j}.
\end{eqnarray}
The last term obviously drops out when $X$ is a symplectic vector field.

If $X$ is a symplectic vector field preserving the connection and
preserving the series of $2$-forms $\Omega$, then $\L_X r = 0$ so
\[
-Di(X)r = i(X)Dr+\textstyle{\frac1\nu}\left[\omega_{ij}X^iy^j +\half
(\nabla_{i}(i(X)
\omega))_{j}y^iy^j +i(X)r~,~r\right]
\]
Using equation (\ref{Dr}), this gives
\[
-Di(X)r = i(X)\ol{R} -i(X)\Omega
+\textstyle{\frac1\nu}\left[\omega_{ij}X^iy^j
+\half (\nabla_{i}(i(X) \omega))_{j}y^iy^j,r\right].
\]
On the other hand, using the fact that
$Da = \partial a - \delta(a) - \frac1{\nu}[r,a]$ one has
\[
D(\omega_{ij}X^iy^j) = -i(X)\omega+ \partial (\omega_{ij}X^iy^j)
+\textstyle{\frac1\nu}[\omega_{ij}X^iy^j,r]
\]
and
\begin{eqnarray*}
D\left(\half (\nabla_{i}(i(X) \omega))_{j}y^iy^j\right)
& = & -\nabla_{i}(i(X) \omega))_{j}dx^iy^j
+ \partial \left(\half (\nabla_{i}(i(X) \omega))_{j}y^iy^j\right)\\
& & \qquad\mbox{}+
\textstyle{\frac1\nu}
\left[\half(\nabla_{i}(i(X) \omega))_{j}y^iy^j,r\right].
\end{eqnarray*}
Since $X$ is an affine vector field, one has
$ (i(X)R) (Y)Z = (\nabla^2 X)(Y,Z)$ so that
\[
\partial \left(\half(\nabla_{i}(i(X) \omega))_{j}y^iy^j\right) =
-\half  ((\nabla^2 X)^p_{ki}\omega)_{jp}y^i y^j dx^k
   = i(X)\ol{R}.
\]
Hence
\[
D\left(-i(X)r - \omega_{ij}X^iy^j - \half (\nabla_{i}(i(X)
\omega))_{j}y^iy^j\right) = i(X) \omega-i(X)\Omega.
\]
So, for any vector field $X$ so that $\L_X\omega = 0,
\L_X \Omega = 0$ and $\L_X \nabla = 0$, one has
\[
\L_X = D\circ i(X) + i(X)\circ D+\textstyle{\frac1\nu}\ad_*(T(X))
\]
with $T(X) = i(X)r + \omega_{ij}X^iy^j +
\half (\nabla_{i}(i(X) \omega))_{j}y^iy^j$ and
\[
DT(X) = -i(X) \omega+i(X)\Omega.
\]
In particular, if there exists a series of smooth
functions $ \lambda_X$
so that
\begin{equation}\label{eq_formalHam}
i(X)\omega-i(X)\Omega = d\lambda_X
\end{equation}
one can write
\[
\L_X = D\circ i(X)+ i(X)\circ D +
\textstyle{\frac1\nu}\ad_*(\lambda_X + T(X))
\]
with
\[
D(\lambda_X+T(X)) = 0.
\]
Thus $\lambda_X+T(X)$ is the flat section associated to the
series of smooth function on $M$ obtained by taking the
part of $\lambda_X+T(X)$ with no $y$ terms hence $\lambda_X$
(notice that $i(X)r$ has
no terms without a $y$ from the construction of $r$).
If $Q$ denotes the quantisation map associating a flat section
to a series in $\nu$ of smooth functions, the above yields
\[
\L_X = D\circ i(X) + i(X)\circ D +
\textstyle{\frac1\nu}\ad_*(Q(\lambda_X)).
\]
Since in those assumptions the map $Q$ commutes with $\L_X$
one has
\[
Q(Xf) = \L_X Q(f) = \textstyle{\frac1\nu}[Q(\lambda_X), Q(f)]
\]
so that for any smooth function $f$, one has
\[
Xf = \textstyle{\frac1\nu}(\ad_* \lambda_X )(f).
\]
This proves Proposition 4.3 of \cite{bib:Kravchenko}.

We now aim to show that the condition (\ref{eq_formalHam}) is not only
sufficient, but also necessary. Observe that any Fedosov star product
has the Poisson bracket for the term of order 1 in $\nu$ and has a
second term which is of order at most $2$ in each argument so, as was
mentioned before, it uniquely defines a symplectic connection (which is
the connection used in the construction) so that invariance of $\nabla$
is a necessary condition for the invariance of $*_{\nabla,\Omega}$. In
\cite{bib:BertCahGutt} it is shown that $\Omega$ can also be recovered,
so in fact we have the following well known Lemma:

\begin{lemma}
A vector field $X$ is a derivation of $*_{\nabla,\Omega}$ if and only if
$\L_X \omega=0$, $\L_X \Omega=0$, and $\L_X \nabla=0$.
\end{lemma}

We have seen above that such a vector field $X$ is an
inner derivation if $i(X) (\omega-\Omega)$ is exact. We shall show now
that this is also a necessary condition.

Assume $X$ is a vector field on $M$ such that
there exists a series of smooth functions $ \lambda_X$ with
\begin{equation}\label{eq_innder}
X(u)=\textstyle\textstyle{\frac1{\nu}}(\ad_*\lambda_X )(u)
\end{equation}
for every smooth function $u$ on $M$.
Then $X$ is a derivation of $*$
so $\L_X \omega=0$, $\L_X \Omega=0$, $\L_X \nabla=0$ and
\[
Q(Xf)= \L_X Q(f)=\textstyle{\frac1{\nu}}[T(X), Q(f)]
\]
with $T(X)=i(X)r+\omega_{ij}X^iy^j+
{\half }(\nabla_{i}(i(X) \omega))_{j}y^iy^j$ and
\[
DT(X)=-i(X) \omega+i(X)\Omega.
\]
Taking a contractible open set $U$ in $M$, there exists a series of
smooth locally defined functions $\lambda^U_X$ on $U$ so that
\[
(i(X)\omega-i(X)\Omega)\vert_U=d\lambda^U_X
\]
and, everything being local, we have on $U$
\[
D(\lambda^U_X+T(X))\vert_U=0,
\]
thus $\lambda^U_X+T(X)$ is the flat section on $U$ associated to the
series of smooth functions on $U$ obtained by taking the part of
$\lambda^U_X+T(X)$ with no $y$ terms (which is $\lambda^U_X$) and
\[
Q(X(u))\vert_U= \L_X Q(u)\vert_U =
\left.\textstyle{\frac1{\nu}}[Q(\lambda^U_X), Q(u)]\right\vert_U
\]
so that
\[
X(u)\vert_U = \left.
\textstyle{\frac1{\nu}}(\ad_{*_{\nabla,\Omega}}
\lambda^U_X )(u)\right\vert_U
\]
for any smooth function $u$. Comparing this with equation
(\ref{eq_innder}) shows that
\[
\lambda^U_X -\lambda_X
\]
is a constant on $U$ and hence that
\[
i(X)\omega-i(X)\Omega=d\lambda_X.
\]
Thus we have proved the converse of Kravchenko's result. In summary:
\begin{theorem}
A vector field $X$ is an inner derivation of $* = *_{\nabla,\Omega}$
if and only if $\L_X \nabla=0$ and there exists a series of functions
$\lambda_X$ such that
\[
i(X)\omega-i(X)\Omega = d \lambda_X.
\]
In this case
\[
X(u)=\textstyle{\frac1\nu}(\ad_*\lambda_X )(u).
\]
\end{theorem}

\section{Moment Maps for an invariant Star Product with an
invariant connection}

Let $(M,\omega)$ be endowed with a differential star product $*$,
\[
u*v =
uv + \sum_{r\ge1} \nu^r C_r(u,v).
\]
Consider an algebra $\g$ of vector fields on $M$ consisting of
derivations of $*$
and assume that there is a symplectic connection $\nabla$ which
is invariant under $\g$ (i.e.\ $\L_X \nabla=0~, \forall X \in \g$).
This is of course automatically true if the star product is natural
and invariant.

It was proven in \cite{bib:BBG} that $*$ is equivalent, through an
equivariant equivalence
\[
T=\Id+\sum_{r\ge1} \nu^r T_r
\]
(i.e.\ $\L_X T=0$), to a Fedosov star product built from $\nabla$
and a series of invariant closed $2$-forms $\Omega$ which give a
representative
of the characteristic class of $*$.

Observe that
\[
X(u)=\textstyle{\frac1{\nu}}(\ad_*\mu_X )(u)
\]
for any $X\in\g$ if and only if
\[
X(u)=T\circ X\circ T^{-1}(u)
= T(\textstyle{\frac1{\nu}}(\ad_*T\mu_X )(T^{-1}u))=
\textstyle{\frac1{\nu}}(\ad_{*_{\nabla,\Omega}}T\mu_X )(u).
\]
Hence the Lie algebra $\g$ consists of inner derivations for $*$ if and
only if this is true for the Fedosov star product $ *_{\nabla,\Omega}$
and this true if and only if there exists a series of functions
$\lambda_X$ such that
\[
i(X)\omega-i(X)\Omega=d\lambda_X.
\]
In this case
\[
X(u)=\textstyle{\frac1{\nu}}(\ad_*\mu_X )(u) \quad \mbox{\rm with}\quad
\mu_X=T\lambda_X.
\]
In particular, this yields
\begin{theorem}
Let $G$ be a compact Lie group of symplectomorphisms of $(M,\omega)$ and
$\g$ the corresponding Lie algebra of symplectic vector fields on $M$.
Consider a star product $*$ on $M$ which is invariant under $G$. The Lie
algebra $\g$ consists of inner derivations for $*$ if and only if there
exists a series of functions $\lambda_X$ and a representative
$\textstyle{\frac1{\nu}}(\omega-\Omega)$ of the characteristic class of
$*$ such that
\[
i(X)\omega-i(X)\Omega = d\lambda_X.
\]
\end{theorem}


\end{document}